\documentclass[a4 paper,12pt,bezier]{article}
\usepackage[top=3cm,bottom=3cm,left=2.5cm,right=2.5cm]{geometry}
\usepackage{tikz, color,xcolor}
\usepackage{caption}
\usepackage{amsmath}
\usepackage{graphics}
\usepackage{mathrsfs}
\usepackage{subcaption}
\usepackage{amsmath}
\usepackage{amsfonts,amsthm,amssymb,mathrsfs,bbding}
\usepackage{float}
\usepackage{graphics,epstopdf}
\usepackage{graphics,multicol}
\usepackage{graphicx}
\usepackage{color}
\usepackage{enumerate}
\usepackage{caption}
\usepackage{rotating}
\usepackage{lscape}
\usepackage{longtable}
\allowdisplaybreaks[4]
\usepackage{tabularx}
\usepackage[colorlinks=true,anchorcolor=blue,filecolor=blue,linkcolor=blue,urlcolor=blue,citecolor=blue]{hyperref}
\usepackage{extarrows}
\usepackage{cite}
\usepackage{latexsym,bm}
\usepackage{mathtools}
\usepackage{changes}
\usepackage{lipsum}
\definechangesauthor[name={kk}, color=blue]{kk}
\newtheorem{theorem}{Theorem}
\newtheorem{prob}{Problem}

\newtheorem{lemma}{Lemma}

\newtheorem{claim}{Claim}

\renewcommand\proofname{\bf Proof}

\definecolor{ECUST}{RGB}{0,106,167}
\addtocounter{section}{0}

\makeatletter
\newsavebox\myboxA
\newsavebox\myboxB
\newlength\mylenA
\newcommand*\xoverline[2][0.75]{%
\sbox{\myboxA}{$\m@th#2$}%
\setbox\myboxB\null
\ht\myboxB=\ht\myboxA%
\dp\myboxB=\dp\myboxA%
\wd\myboxB=#1\wd\myboxA
\sbox\myboxB{$\m@th\overline{\copy\myboxB}$}
\setlength\mylenA{\the\wd\myboxA}
\addtolength\mylenA{-\the\wd\myboxB}%
\ifdim\wd\myboxB<\wd\myboxA%
\rlap{\hskip 0.5\mylenA\usebox\myboxB}{\usebox\myboxA}%
\else
\hskip -0.5\mylenA\rlap{\usebox\myboxA}{\hskip 0.5\mylenA\usebox\myboxB}%
\fi}
\makeatother
\begin{document}
\title{\bf Spectral Tur\'{a}n  problem for $\mathcal{K}_5^{-}$-free signed graphs}
\author{Yongang Wang\footnote{Email: yongang\_wang@163.com}\\[2mm]
\small School of Mathematics, East China University of Science and Technology, \\
\small  Shanghai 200237, P.R. China\\}
\date{}
\maketitle
{\flushleft\large\bf Abstract}

The classical spectral Tur\'{a}n  problem is to determine the maximum spectral radius of an  $\mathcal{F}$-free graph of order $n.$ Let  $\mathcal{K}_{k}^{-}$ be the set of all   unbalanced $K_k.$ In this paper, we focus on the spectral Tur\'{a}n  problem of $\mathcal{K}_{k}^{-}$-free unbalanced signed graph for $k\geq5$. Moreover, we give an answer for $k=5$  and completely characterize the corresponding extremal signed graph.

\begin{flushleft}
\end{flushleft}

\section{Introduction}

 All graphs in this paper are simple. A  \textit{signed graph}   $\Gamma=(G,\sigma)$  consists of a graph $G=(V(G),E(G))$ and a sign function $\sigma : E\rightarrow \{-1,1\}$, where $G$ is its underlying graph and  $\sigma$ is its sign function. An edge $e$ is \textit{positive}\ (\textit{negative})\ if $\sigma(e)=1$ (resp.  $\sigma(e)=-1$).  A cycle  $C$ in  a signed graph $\Gamma$ is called \emph{positive} (resp.\ \emph{negative})  if the number of its negative edges is even (resp.\ odd). A signed graph is called \emph{balanced} if  all  its cycles are positive; otherwise, it is called \emph{unbalanced}.  The \textit{adjacency} \textit{matrix} of  $\Gamma$ is denoted by $A(\Gamma)=(a^{\sigma}_{ij})$,\ where $a^{\sigma}_{ij} =\sigma(v_{i}v_{j})$ if $v_{i}\sim v_{j}$,\ and $0$ otherwise. The eigenvalues of $A({\Gamma})$ are called the eigenvalues of  $\Gamma,$ denoted by $\lambda_1(\Gamma)\geq \lambda_2(\Gamma)\geq \ldots\geq \lambda_n(\Gamma).$   In particular,   $\lambda_{1}(\Gamma)$ is called the \textit{index} of $\Gamma.$ The \textit{spectral radius} of $\Gamma$ is defined by $\rho(\Gamma)=\max\{|\lambda_i(\Gamma)|:1\leq i\leq n\}.$ 
 For more details about the notion of signed graphs, we refer to\cite{treepositive,cycle}.

The   spectral theory of signed
 graph  has been studied extensively in the literature. For the largest eigenvalue  of a signed graph with   certain structures, 
 Koledin and Stani\'{c}\cite{connectedindex} studied connected signed graphs of fixed order, size and number of negative edges that maximize the index  of their adjacency matrices. After that, signed graphs maximizing the index in suitable subsets of  signed complete graphs have been studied by Ghorbani and Majidi\cite{maxindex},  Li, Lin and Meng\cite{Meng} and   Akbari,  Dalvandi,  Heydari and  Maghasedi\cite{kedge}.   It is well known that the eigenvalues of a balanced signed graph are the same as those of its underlying graph. Therefore, the largest eigenvalue  of an unbalanced signed graph  has attracted more attention of scholars.   In  $2019$, Akbari, Belardo, Heydari, Maghasedi and Souri\cite{unicyclic}  determined the signed graphs achieving the minimal or  maximal index in the class of unbalanced signed unicyclic graphs. In 2021,   He, Li, Shan and Wang \cite{shuangquan} gave the first five largest indices among all unbalanced signed bicyclic graphs of order  $n\geq36$.    In $2022$,  Brunetti and Stani\'{c}\cite{unbalancedindex} studied the extremal spectral radius among  all unbalanced connected signed graphs.   More results on the  the spectral theory of signed graphs can be found in \cite{open,signedturan,WYQ,JCTA,huanghao}, where \cite{open} is an excellent survey about some general results and  problems on the spectra of signed graphs.

 Let $\mathcal{F}$ be a family of graphs. A graph $G$ is $\mathcal{F}$-free if $G$ does not contain any graph in $\mathcal{F}$
as a subgraph. The classical spectral Tur\'{a}n  problem is to determine the maximum spectral radius of an  $\mathcal{F}$-free graph of order $n$, which is known as  the \textit{spectral Tur\'{a}n number} of $\mathcal{F}$.  This problem was originally proposed by Nikiforov\cite{proposeNikiforov}.  With regard to  unsigned graphs,  Nikiforov \cite{NikiKr} gave the  spectral   Tur\'{a}n number of  complete graph $K_r.$  In the wake of that,  much attention has been paid to the spectral Tur\'{a}n  problem in the past decades, such as $K_{s,t}$-free graphs \cite{1Babai,NikiKr},  $C_4$-free graphs \cite{NikiC4,NikiKr,zhaiC4},  $C_6$-free graphs \cite{C6Zhai} and $C_{2k+1}$-free graphs \cite{oddNiki}.   For more excellent results, we refer the reader to \cite{Li,triangleLin,Bollobas,WilfKr,JCTBlin}. In this paper, we focus on the spectral Tur\'{a}n  problem in signed graphs.


The spectral Tur\'{a}n problem in signed graphs  has been studied in recent years.  Wang and Lin \cite{my}  gave  a spectral  condition for the existence of negative $C_4$ in  unbalanced signed graphs.  Let   $\mathcal{K}_k^{-}$ be the set of  all unbalanced $K_k$. In 2022,  Wang, Hou and Li \cite{C3free} determined the spectral Tur\'{a}n number of $\mathcal{K}_3^{-}.$  In 2023,      Chen and Yuan\cite{K4free} gave the spectral Tur\'{a}n number of $\mathcal{K}_4^{-}.$  Motivated by their work,  we focus on studying the existence of unbalanced $K_{k}$ for $k\geq 5$ (see Problem \ref{problem}). 

\begin{prob}\label{problem}
What is the maximum spectral radius among all $\mathcal{K}_{k}^{-}$-free unbalanced signed graphs   for $k\geq 5?$
\end{prob}

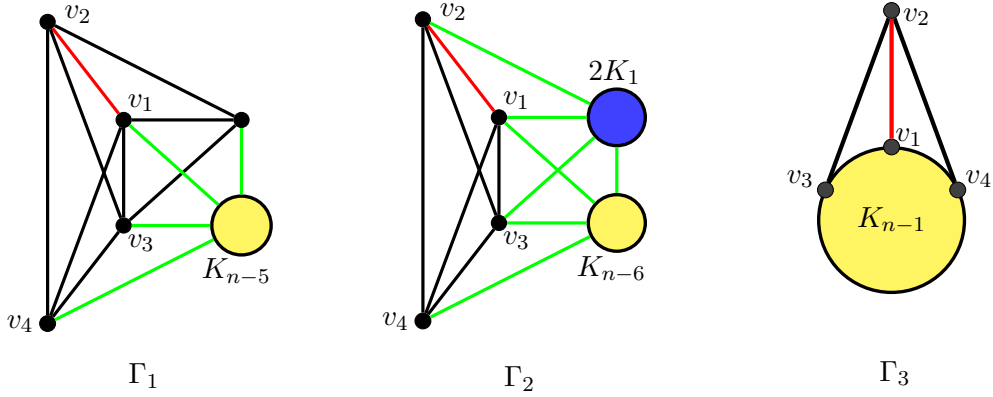
\begin{figure}[htbp!]
  \centering

\begin{subfigure}[c]{0.3\linewidth}
    \centering
    \begin{tikzpicture}
      \tikzstyle{every node}=[font=\small]

       \node (1) at (0,0) [shape=circle,draw,fill=black!100,inner sep=2.1pt,minimum size=2mm] {};
      \node (2) at (0,4) [shape=circle,draw,fill=black!100,inner sep=1pt,minimum size=2mm] {};
      \node (3) at (1,1.3) [shape=circle,draw,fill=black!100,inner sep=1pt,minimum size=2mm] {};
      \node (4) at (1,2.7) [shape=circle,draw,fill=black!100,inner sep=1pt,minimum size=2mm] {};
      \node (7) at (0.4,4.1) [] {$v_2$};
       \node (8) at (-0.36,0) [] {$v_4$};
        \node (9) at (1.23,2.93) [] {$v_1$};
         \node (10) at (1.23,1.07) [] {$v_3$};
          \node (11) at (2.48,0.69) [] {$K_{n-5}$};
     \node (5) at (2.55,1.3) [shape=circle,line width=1.26pt,draw,fill=yellow!75,inner sep=1pt,minimum size=7.8mm] {};
     \node (6) at (2.55,2.7) [shape=circle,draw,fill=black!100,inner sep=1pt,minimum size=2mm] {};
      \draw [red,line width=1.2pt] (4) to (2);
       \draw [black,line width=1.2pt] (3) to (1);
      \draw [black,line width=1.2pt] (3) to (2);
            \draw [black,line width=1.2pt] (4) to (1);
      \draw [black,line width=1.2pt] (1) to (2);
            \draw [black,line width=1.2pt] (3) to (4);
      \draw [black,line width=1.2pt] (6) to (2);
            \draw [green,line width=1.2pt] (5) to (1);
      \draw [green,line width=1.2pt] (3) to (5);
            \draw [black,line width=1.2pt] (3) to (6);
               \draw [green,line width=1.2pt] (4) to (5);
            \draw [black,line width=1.2pt] (4) to (6);
      \draw [green,line width=1.2pt] (5) to (6);
  \end{tikzpicture}
\captionsetup{labelformat=empty}
\vspace{-1pt}
    \subcaption{$\Gamma_1$}
\end{subfigure}
\begin{subfigure}[c]{0.3\linewidth}
    \centering
    \begin{tikzpicture}
      \tikzstyle{every node}=[font=\small]
       \node (1) at (0,0) [shape=circle,draw,fill=black!100,inner sep=2.1pt,minimum size=2mm] {};
      \node (2) at (0,4) [shape=circle,draw,fill=black!100,inner sep=1pt,minimum size=2mm] {};
      \node (3) at (1,1.3) [shape=circle,draw,fill=black!100,inner sep=1pt,minimum size=2mm] {};
      \node (4) at (1,2.7) [shape=circle,draw,fill=black!100,inner sep=1pt,minimum size=2mm] {};
      \node (7) at (0.4,4.1) [] {$v_2$};
       \node (8) at (-0.36,0) [] {$v_4$};
        \node (9) at (1.23,2.93) [] {$v_1$};
         \node (10) at (1.23,1.07) [] {$v_3$};
          \node (11) at (2.48,0.69) [] {$K_{n-6}$};
           \node (12) at (2.52,3.3) [] {$2K_1$};
     \node (5) at (2.55,1.3) [shape=circle,line width=1.26pt,draw,fill=yellow!75,inner sep=1pt,minimum size=7.5mm] {};
     \node (6) at (2.55,2.7)  [shape=circle,line width=1.26pt,draw,fill=blue!75,inner sep=1pt,minimum size=7.5mm] {};
      \draw [red,line width=1.2pt] (4) to (2);
       \draw [black,line width=1.2pt] (3) to (1);
      \draw [black,line width=1.2pt] (3) to (2);
            \draw [black,line width=1.2pt] (4) to (1);
      \draw [black,line width=1.2pt] (1) to (2);
            \draw [black,line width=1.2pt] (3) to (4);
      \draw [green,line width=1.2pt] (6) to (2);
            \draw [green,line width=1.2pt] (5) to (1);
      \draw [green,line width=1.2pt] (3) to (5);
            \draw [green,line width=1.2pt] (3) to (6);
               \draw [green,line width=1.2pt] (4) to (5);
            \draw [green,line width=1.2pt] (4) to (6);
      \draw [green,line width=1.2pt] (5) to (6);
  \end{tikzpicture}
\captionsetup{labelformat=empty}
\vspace{1pt}
    \subcaption{$\Gamma_2$}
\end{subfigure}
    \begin{subfigure}[c]{0.3\linewidth}
    \centering
    \begin{tikzpicture}
      \tikzstyle{every node}=[font=\small]
 \node (1) at (0,0) [shape=circle,line width=1.26pt,draw,fill=yellow!75,inner sep=10pt,minimum size=10.7mm] {$K_{n-1}$};
      \node (2) at (0.87,0.4) [shape=circle,draw,fill=black!75,inner sep=2.1pt,minimum size=1.2mm] {};
      \node (3) at (-0.87,0.4) [shape=circle,draw,fill=black!75,inner sep=2.1pt,minimum size=1.2mm] {};
      \node (4) at (0,0.97) [shape=circle,draw,fill=black!75,inner sep=2.1pt,minimum size=1.2mm] {};
      \node (5) at (0,2.78) [shape=circle,draw,fill=black!75,inner sep=2.1pt,minimum size=1.2mm] {};
     \node (11) at (0.23,1.1) [] {$v_1$};
      \node (12) at (-1.2,0.55) [] {$v_3$};
      \node (13) at (1.15,0.55) [] {$v_4$};
       \node (13) at (0.33,2.7) [] {$v_2$};
      \draw [black,line width=1.6pt] (5) to (2);
      \draw [black,line width=1.6pt] (5) to (3);
      \draw [red,line width=1.6pt] (5) to (4);
  \end{tikzpicture}
\captionsetup{labelformat=empty}
\vspace{16.3pt}
    \subcaption{$\Gamma_3$}
\end{subfigure}
	\caption{The signed graphs $\Gamma_1,$  $\Gamma_2$ and $\Gamma_3$.}\label{fig123}
\end{figure}

Suppose $\Gamma=(G,\sigma)$ is a signed graph and $U\subset V(G)$.
The operation that changes the signs of all edges between $U$ and $V(G)\backslash U$ is called a \emph{switching operation}. If a signed graph $\Gamma^{\prime}$ is obtained from $\Gamma$ by applying finitely many switching operations, then $\Gamma$ is said to be \emph{switching equivalent} to $\Gamma^{\prime}$.  For $n\geq7$, we define the signed graphs $\Gamma_1,$  $\Gamma_2$ and $\Gamma_3$ as shown in Figure \ref{fig123}, where the yellow and blue circles represent (unsigned) complete graphs and independent set,  respectively, the red lines represent negative edges and the other lines represent positive edges, and especially,  the green lines joined to the circles represent the connection of all possible edges. In this paper, we give an answer to Problem \ref{problem} for $k=5$ as follows.

\begin{theorem}\label{C4-free}
Let $\Gamma=(G,\sigma)$ be an unbalanced signed graph of order $n\geq5$. If $$\rho(\Gamma)\geq\lambda_{1}(\Gamma_{3}),$$  then  $\Gamma$ contains an unbalanced  $K_5$ unless $\Gamma$ is switching equivalent to $\Gamma_{3}$ (see Figure \ref{fig123}).
\end{theorem}

\section{The largest eigenvalue of signed graph $\Gamma_3$}
In this section, we shall show that $\lambda_{1}(\Gamma_3)>n-2.$  We now introduce the definition of equitable  quotient matrix.

Let $M$ be a real symmetric matrix of order $n$, and let $[n]=\{1,2, \ldots, n\}$. Given a partition $\Pi:[n]=X_1 \cup X_2 \cup \cdots \cup X_k$, the matrix $M$ can be written as
$$
M=\left[\begin{array}{cccc}
M_{1,1} & M_{1,2} & \cdots & M_{1, k} \\
M_{2,1} & M_{2,2} & \cdots & M_{2, k} \\
\vdots & \vdots & \ddots & \vdots \\
M_{k, 1} & M_{k, 2} & \cdots & M_{k, k}
\end{array}\right] .
$$
If all row sums of $M_{i, j}$ are the same, say $b_{i, j}$, for all $i, j \in\{1,2, \ldots, k\}$, then $\Pi$ is called an $equitable$   $partition$ of $M$, and the matrix $Q=\left(b_{i, j}\right)_{i, j=1}^k$ is called an $equitable$  $quotient$  $matrix$ of $M$.

\begin{lemma}{\bf(\cite[p.24]{fz})}\label{quotient}
Let $M$ be a real symmetric matrix, and let $Q$ be an equitable quotient matrix of $M$.  Then the matrix $M$
 has the following two kinds of eigenvalues.
 
 $(1)$ The eigenvalues coincide with the eigenvalues of $Q.$
 
 $(2)$ The eigenvalues of $M$ not in the spectra of $Q$ remain unchanged if some scalar multiple of
 the all-one block J is added to block $M_{ij}$ for each $1 \leq i,j \leq k.$
\end{lemma}

\begin{lemma}\label{dayu n3}
Let $\Gamma_1,$  $\Gamma_2$  and $\Gamma_3$ be the signed graphs as shown in Figure \ref{fig123}. Then we have the following statements.
\begin{itemize}
\item[(1)] $\lambda_{1}(\Gamma_1)< n-2.$
\item[(2)] $\lambda_{1}(\Gamma_2)< n-2.$
\item[(3)] $\lambda_{1}(\Gamma_3)>n-2.$
\end{itemize}
\end{lemma}
\renewcommand\proofname{\bf Proof}
\begin{proof}
$(1)$ Let  $J$, $I$  and $O$ denote the  all-ones matrix, identity matrix and all-zeros matrix, respectively.  By a suitable partition,
$$A(\Gamma_1)=\left[\begin{array}{ccccc}
     0 & -1 & 1 & J & J\\
     -1 & 0 & 1 & J & O\\
      1 & 1 & 0 & J & J\\
     J & J &  J & O_2 & J\\
     J & O & J & J & (J-I)_{n-5}\\
\end{array}\right],$$
and $A(\Gamma_1)$ has the equitable quotient matrix 
$$Q_1=\left[\begin{array}{ccccc}
       0 & -1 & 1 & 2 & n-5\\
     -1 & 0 & 1 & 2 & 0\\
      1 & 1 & 0 & 2 & n-5\\
     1 & 1 &  1 & 0 & n-5\\
     1 & 0 & 1 & 2 & n-6\\
\end{array}\right].$$
By a simple calculation, the characteristic polynomial of $Q_1$ is
$$g(x)=x^5+(6-n)x^4+(11-4n)x^3-(6+n)x^2+(8n-32)x+6n-20.$$ Observe that  $g(-\infty)<0,$ $g(-2)>0,$  $g(-\frac{6}{5})<0,$  $g(0)>0,$  $g(n-3)<0$ and  $g(n-2)>0.$ Then  $\lambda_1(\Gamma_1)<n-2.$  By adding  some scalar multiple of the all-one block $J$ to the block of $A(\Gamma_1),$  $A(\Gamma_1)$  becomes
$$A_1=\left[\begin{array}{ccccc}
     0 & 0 & 0 & 0 & 0\\
     0 & 0 & 0 & 0 & O\\
      0 & 0 & 0 & 0 & 0\\
     0 & 0 &  0 & O_2 & 0\\
     0& O & 0 & 0 & -I_{n-5}\\
\end{array}\right].$$  Note that the distinct eigenvalues of $A_1$ are $-1$ and $0.$  Then by Lemma \ref{quotient}, $\lambda_1(\Gamma_1)=\lambda_1(Q_1)<n-2.$

$(2)$ By a suitable partition,
$$A(\Gamma_2)=\left[\begin{array}{ccccc}
     0 & -1 & 1 & J & J\\
     -1 & 0 & 1 & J & O\\
      1 & 1 & 0 & J & J\\
     J & J &  J & O_3 & J\\
     J & O & J & J & (J-I)_{n-6}\\
\end{array}\right],$$
and $A(\Gamma_2)$ has the equitable quotient matrix 
$$Q_2=\left[\begin{array}{ccccc}
       0 & -1 & 1 & 3 & n-6\\
     -1 & 0 & 1 & 3 & 0\\
      1 & 1 & 0 & 3 & n-6\\
     1 & 1 &  1 & 0 & n-6\\
     1 & 0 & 1 & 3 & n-7\\
\end{array}\right].$$
By a simple calculation, the characteristic polynomial of $Q_2$ is
$$h(x)=x^5+(7-n)x^4+(18-5n)x^3-(4+2n)x^2+(11n-55)x+9n-39.$$  Observe that  $h(-\infty)<0,$ $h(-2)>0,$  $h(-\frac{6}{5})<0,$  $h(0)>0,$  $h(n-3)<0$ and  $h(n-2)>0.$ Then  $\lambda_1(\Gamma_2)<n-2.$  By adding  some scalar multiple of the all-one block $J$ to the block of $A(\Gamma_2),$  $A(\Gamma_2)$  becomes
$$A_2=\left[\begin{array}{ccccc}
     0 & 0 & 0 & 0 & 0\\
     0 & 0 & 0 & 0 & O\\
      0 & 0 & 0 & 0 & 0\\
     0 & 0 &  0 & O_3 & 0\\
     0& O & 0 & 0 & -I_{n-6}\\
\end{array}\right].$$  Note that the distinct eigenvalues of $A_2$ are $-1$ and $0.$  Then by Lemma \ref{quotient}, $\lambda_1(\Gamma_2)=\lambda_1(Q_2)<n-2.$ 

$(3)$ By a suitable partition,
$$A(\Gamma_3)=\left[\begin{array}{cccc}
     0 & -1 & J & O\\
     -1 & 0& J & J\\
     J & J& (J-I)_2 & J\\
     O & J& J& J-I\\
\end{array}\right],$$
and $A(\Gamma_1)$ has the equitable quotient matrix
$$Q_3=\left[\begin{array}{cccc}
     0 & -1 & 2 & 0\\
     -1 & 0& 2  & n-4\\
     1 & 1& 1 & n-4\\
     0 & 1& 2& n-5\\
\end{array}\right].$$
By a simple calculation, the characteristic polynomial of $Q_3$ is
$$f(x)=(x+1)(x^3+(3-n)x^2-(n+1)x+3n-7).$$
 Observe that $f(n-2)=1-n<0.$ Then $\lambda_1(\Gamma_3)\geq\lambda_1(Q_3)>n-2$ by Lemma \ref{quotient}.

\end{proof}

\section{Proof of Theorem \ref{C4-free}}
By computer searching and the table of the spectra of signed graphs with at most six vertices\cite{table}, we can check that Theorem \ref{C4-free} is true for $n\leq 8.$ Therefore, we now assume that $n\geq 9.$  We first give two lemmas  which are needed  in the proof of Theorem \ref{C4-free}. 

\begin{lemma}{\bf (\cite[Lemma 2.5]{SGX})}\label{Nonnegative}
Let $\Gamma$ be a signed graph. Then there exists a signed graph $\Gamma^{\prime}$ switching equivalent to $\Gamma$ such that   $A(\Gamma^{\prime})$ has a non-negative eigenvector corresponding to $\lambda_{1}(\Gamma^{\prime})$.
\end{lemma}

\begin{lemma}{\bf (\cite[Proposition 3.2]{treepositive})}\label{remaincycle}
 Two signed graphs with the same underlying graph are switching equivalent
 if and only if they have the same set of positive cycles.
\end{lemma}

 Wang,  Hou and  Li and\cite{C3free} and Chen and Yuan\cite{K4free} determined the spectral Tur\'{a}n number of $\mathcal{K}_3^{-}$ and $\mathcal{K}_4^{-}$ in unbalanced signed graphs, respectively.

\begin{lemma}{\bf (\cite[Theorem 1.3]{C3free})}\label{C3-free}
   Let $\Gamma=(G,\sigma)$ be a connected unbalanced signed graph of order $n$. If $\Gamma$ is $\mathcal{K}_{3}^{-}$-free, then $$\rho(\Gamma)\leq\frac{1}{2}(\sqrt{n^2-8}+n-4).$$
\end{lemma}

\begin{lemma}{\bf (\cite[Theorem 1.6]{K4free})}\label{K4-free}
   Let $\Gamma=(G,\sigma)$ be an unbalanced signed graph of order $n$. If $\Gamma$ is $\mathcal{K}_{4}^{-}$-free, then $$\rho(\Gamma)\leq n-2.$$
\end{lemma}

For a signed graph $\Gamma$, the frustration index, denoted by $ \epsilon(\Gamma)$, is the minimum number of edges to be deleted such that the resultant signed graph is balanced\cite{frustration}.
\begin{lemma}{\bf (\cite[Theorem 3.2]{2m-n+1})}\label{2n-n+1}
Let $\Gamma$ be a connected signed graph with $n$ vertices, $m$ edges and the frustration index $\epsilon(\Gamma).$ Then $$\lambda_1(\Gamma)\leq \sqrt{2(m-\epsilon(\Gamma))-n+1}.$$
\end{lemma}
The  clique number of a graph $G$, denoted by $\omega(\Gamma)$, is the maximum order of a  clique in $G.$ The balanced clique number of a signed graph $\Gamma$, denoted by $\omega_b(\Gamma)$, is the maximum order of a balanced clique in $\Gamma.$

\begin{lemma}{\bf (\cite{WilfKr})}\label{Gseshu}
Let $G$ be a  graph of order $n$. Then
$$\lambda_1(G)\leq n \left (1-\frac{1}{\omega(G)}\right).$$
\end{lemma}

\begin{lemma}{\bf (\cite[Proposition 5]{WYQ})}\label{seshu}
Let $\Gamma$ be a signed graph of order $n$. Then
$$\lambda_1(\Gamma)\leq n \left (1-\frac{1}{\omega_b(\Gamma)}\right).$$
\end{lemma}
Now, we  are in a position to give the proof of Theorem \ref{C4-free}.
\renewcommand\proofname{\bf Proof of Theorem \ref{C4-free}}
\begin{proof}
Suppose that $\Gamma=(G,\sigma)$ has the maximum spectral radius among all   $\mathcal{K}_{5}^{-}$-free unbalanced signed graphs. We  shall show that $\Gamma$ is switching equivalent to $\Gamma_3.$  Recall that $\rho(\Gamma)=\max\{\lambda_{1}(\Gamma),-\lambda_{n}(\Gamma)\}.$  Since $\Gamma_{3}$ is unbalanced and $\mathcal{K}_{5}^{-}$-free, we may suppose that  $\rho(\Gamma)\geq\lambda_1(\Gamma_{3})>n-2$ by Lemma \ref{dayu n3}.   First we give  some claims. 
\begin{claim}\label{spectral111}
$\rho(\Gamma)=\lambda_{1}(\Gamma).$ 
\end{claim}
Otherwise, assume $\rho(\Gamma)=-\lambda_{n}(\Gamma).$ Let $\Gamma_4=-\Gamma.$ Then $\lambda_{1}(\Gamma_4)=-\lambda_{n}(\Gamma).$  Since $\Gamma$ is  $\mathcal{K}_{5}^{-}$-free, we have   $\omega_{b}(\Gamma_4)\leq4.$  By Lemma \ref{seshu},   for $n\geq 9,$
$$\rho(\Gamma)=-\lambda_{n}(\Gamma)=\lambda_{1}(\Gamma_4)\leq n \left (1-\frac{1}{\omega_b(\Gamma_4)}\right)\leq \frac{3}{4}n<n-2,$$ a contradiction. So Claim \ref{spectral111} holds.


Let  $\Gamma^{\prime}=(G,\sigma^{\prime})$ be a signed graph switching equivalent to $\Gamma.$ By Lemma \ref{Nonnegative}, we can assume that  $A(\Gamma^{\prime})$ has a non-negative eigenvector corresponding to $\lambda_{1}(\Gamma^{\prime})$.   Then by Lemma \ref{remaincycle},     $\Gamma^{\prime}$ is also unbalanced and $\mathcal{K}_{5}^{-}$-free. Furthermore,  $\Gamma^{\prime}$ also has the maximum spectral radius among all $\mathcal{K}_{5}^{-}$-free unbalanced  signed graphs, and  $ \lambda_1(\Gamma^{\prime})=\rho(\Gamma^{\prime})=\rho(\Gamma)>n-2.$  Set $V(\Gamma^{\prime})=\{v_1,v_2,\ldots,v_n\}.$
Let $x=(x_1,x_2,\ldots,x_n)^{\top}$ be the non-negative unit eigenvector of $A(\Gamma^{\prime})$ corresponding to  $\lambda_1(\Gamma^{\prime}),$ where $x_i$ corresponds to the vertex $v_i$ for $1\leq i\leq n.$  Then $$\lambda_1(\Gamma^{\prime})=x^{\top}A(\Gamma^{\prime})x.$$
Now, we begin to analyze the structure of $\Gamma^{\prime}$. 
\begin{claim}\label{1zero}\label{Claim1}
$x$ is a positive vector.
\end{claim}
Otherwise, without loss of generality, assume $x_1=0,$ then $$\begin{aligned}
\lambda_1(\Gamma^{\prime})&=x^{\top}A(\Gamma^{\prime})x=(x_2,x_3,\ldots,x_n)A(\Gamma^{\prime}-v_1)(x_2,x_3,\ldots,x_n)^{\top}\\
&\leq\lambda_1(\Gamma^{\prime}-v_1)\leq\lambda_1(K_{n-1})=n-2,
\end{aligned} $$
a contradiction. So Claim \ref{Claim1} holds.

\begin{claim}\label{connected}
  $\Gamma^{\prime}$ is connected.
\end{claim}
Otherwise, assume $\Gamma^{\prime}_1$ and $\Gamma^{\prime}_2$ are two distinct connected components of $\Gamma^{\prime}$, where $\lambda_1(\Gamma^{\prime})=\lambda_1(\Gamma^{\prime}_1)$.  Without loss of generality, we choose two vertices $v_i\in V(\Gamma^{\prime}_1)$ and $v_j\in V({\Gamma^{\prime}_2}).$  Then we can construct a new signed graph $\Gamma^{\ast}$   obtained from $\Gamma^{\prime}$ by adding a positive edge $v_iv_j.$  Clearly, $\Gamma^{\ast}$ is also unbalanced and $\mathcal{K}_5^-$-free. By Rayleigh principle and Claim \ref{1zero}, we obtain that
$$\begin{aligned}
\lambda_1(\Gamma^{\ast})-\lambda_1(\Gamma^{\prime})&\geq x^{\top}A(\Gamma^{\ast})x-x^{\top}A(\Gamma^{\prime})x\\
&=2x_ix_j>0.
\end{aligned} $$
Hence, $\lambda_1(\Gamma^{\ast})>\lambda_1(\Gamma^{\prime}),$ a contradiction.

Since $\Gamma^{\prime}$ is unbalanced,   $\Gamma^{\prime}$ contains at least one  negative edge and at least one negative cycle. Let $\mathscr{C}$ be one of the shortest negative cycles of $\Gamma^{\prime}$.

\begin{claim}\label{negative triangle}
$\mathscr{C}$ is a negative triangle.
\end{claim}
Otherwise, assume that  $\Gamma^{\prime}$ is $\mathcal{K}_3^-$-free, then by Lemma \ref{C3-free},$$\lambda_{1}(\Gamma^{\prime})\leq\frac{1}{2}(\sqrt{n^2-8}+n-4)<n-2,$$  a contradiction. So Claim \ref{negative triangle} holds.

\begin{claim}\label{containallnegative}
$\mathscr{C}$ contains all negative edges of  $\Gamma^{\prime}.$
\end{claim}
Otherwise,  without loss of generality,  assume $e=v_iv_j$ is a negative edge of $\Gamma^{\prime}$ and $e\notin E(\mathscr{C}).$ Then we can construct a new signed graph $\Gamma^{\ast}$  obtained from $\Gamma^{\prime}$ by deleting  $e.$  Clearly, $\Gamma^{\ast}$ is also unbalanced and $\mathcal{K}_5^-$-free. By Rayleigh principle and Claim \ref{1zero}, we obtain that
$$\begin{aligned}
\lambda_1(\Gamma^{\ast})-\lambda_1(\Gamma^{\prime})&\geq x^{\top}A(\Gamma^{\ast})x-x^{\top}A(\Gamma^{\prime})x\\
&=2x_ix_j>0.
\end{aligned} $$
  Hence, $\lambda_1(\Gamma^{\ast})>\lambda_1(\Gamma^{\prime}),$ a contradiction.
  \begin{claim}\label{containK5}
$G$ contains $K_5$ as a subgraph. 
\end{claim}
  Otherwise, assume $G$ is $K_5$-free, then $\omega(G)\leq4.$ Thus, by Lemma \ref{Gseshu}, 
  $$\lambda_1(\Gamma)\leq \lambda_1(G)\leq n \left (1-\frac{1}{\omega(G)}\right)=\frac{3}{4}n <n-2,$$ a contradiction.

 \begin{claim}\label{ALLedgeK5positive}
The edges of any signed $K_5$ in $\Gamma^{\prime}$ are all positive.
\end{claim}
Otherwise,   assume that $K_5^{\prime}$ is a signed $K_5$ in $\Gamma^{\prime}$ and   $K_5^{\prime}$ contains at least one  negative edge. By Claims \ref{negative triangle}  and \ref{containallnegative},  $K_5^{\prime}$ contains at most three negative edges. If $K_5^{\prime}$ contains three negative edges, then $K_5^{\prime}$ contains $\mathscr{C}$  as a signed subgraph, and thus $K_5^{\prime}$ is unbalanced, a contradiction. If $K_5^{\prime}$ contains one or two negative edges, then $K_5^{\prime}$ must contain a triangle with one negative edge, and thus $K_5^{\prime}$ is unbalanced, a contradiction.

\begin{claim}\label{only1negative}\
$\Gamma^{\prime}$ contains exactly one negative edge.
\end{claim}
Otherwise, assume that $\Gamma^{\prime}$ contains  three negative edges by Claims \ref{negative triangle}  and \ref{containallnegative}.  Without loss of generality, let $v_1v_2$ and $v_2v_3$ be two negative edges of $\Gamma^{\prime}.$ Then  we can construct a new signed graph $\Gamma^{\ast}$  obtained from $\Gamma^{\prime}$ by reversing the sign of   $v_1v_2$ and $v_2v_3.$  By Claim \ref{ALLedgeK5positive}, $\Gamma^{\ast}$ is also unbalanced and $\mathcal{K}_5^-$-free.  Furthermore, by Rayleigh principle and Claim \ref{1zero}, we obtain that
$$\begin{aligned}
\lambda_1(\Gamma^{\ast})-\lambda_1(\Gamma^{\prime})&\geq x^{\top}A(\Gamma^{\ast})x-x^{\top}A(\Gamma^{\prime})x\\
&=2(x_1x_2+x_2x_3)-2(-x_1x_2-x_2x_3)\\
&=4x_2(x_1+x_3)>0.
\end{aligned} $$
Hence, $\lambda_1(\Gamma^{\ast})>\lambda_1(\Gamma^{\prime}),$ a contradiction.

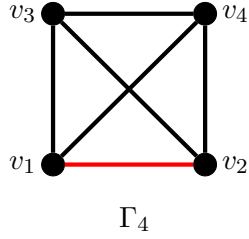
\begin{figure}[htbp!]
  \centering
\begin{subfigure}[c]{0.35\linewidth}
    \centering
    \begin{tikzpicture}
      \tikzstyle{every node}=[font=\small]
  \node (1) at (0,0) [shape=circle,draw,fill=black!100,inner sep=2.1pt,minimum size=3mm] {};
      \node (2) at (2,0) [shape=circle,draw,fill=black!100,inner sep=2.1pt,minimum size=3mm] {};
      \node (3) at (0,2) [shape=circle,draw,fill=black!100,inner sep=2.1pt,minimum size=3mm] {};
      \node (4) at (2,2) [shape=circle,draw,fill=black!100,inner sep=2.1pt,minimum size=3mm] {};
     \node (11) at (-0.4,0) [] {$v_1$};
      \node (12) at (2.4,0) [] {$v_2$};
      \node (13) at (-0.4,2) [] {$v_3$};
       \node (14) at (2.4,2) [] {$v_4$};
      \draw [red,line width=1.6pt] (1) to (2);
      \draw [black,line width=1.6pt] (2) to (4);
      \draw [black,line width=1.6pt] (1) to (4);
        \draw [black,line width=1.6pt] (3) to (4);
      \draw [black,line width=1.6pt] (1) to (3);
        \draw [black,line width=1.6pt] (2) to (3);
  \end{tikzpicture}
\captionsetup{labelformat=empty}
\vspace{0pt}
    \subcaption{$\Gamma_4$}
\end{subfigure}
	\caption{The signed graph $\Gamma_5$.}\label{fig45}
\end{figure}

 By Claim \ref{only1negative},   $\epsilon(\Gamma^{\prime})=1.$ Without loss of generality,  we can suppose that $v_1v_2$ is the unique negative edge of $\Gamma^{\prime}.$  Let $m=|E(\Gamma^{\prime})|=\frac{n(n-1)}{2}-h.$  
\begin{claim}\label{kfreebianshu}
$h\leq n-3.$
\end{claim} 
Otherwise, assume 
$h\geq n-2.$ Then by Lemma \ref{2n-n+1},
$$\begin{aligned}
\lambda_1(\Gamma^{\prime})
&\leq \sqrt{2(m-\epsilon(\Gamma^{\prime}))-n+1}\\
& = \sqrt{2\left(\left(\frac{n(n-1)}{2}-h\right)-1\right)-n+1}\\
&\leq \sqrt{2\left(\left(\frac{n(n-1)}{2}-(n-2)\right)-1\right)-n+1}\\
&=\sqrt{n(n-1)-2n+4-2-n+1}\\
&=\sqrt{n^2-4n+3}\\
&< n-2,
\end{aligned} $$
a contradiction. So Claim \ref{kfreebianshu} holds.

Now,  we define the  signed graph $\Gamma_{5}$   as shown in Figure \ref{fig45}, where the black and red lines represent positive and negative edges,  respectively. By Lemma \ref{K4-free},  $\Gamma^{\prime}$  contains  unbalanced $K_4$ as a signed subgraph.  Without loss of generality,  let $K_4^{-}$ be one of the unbalanced $K_4$ in $\Gamma$  and $V(K_4^-)=\{v_1,v_2,v_3,v_4\}.$ Then by Claim \ref{only1negative}, $K_4^-=\Gamma_5$.  Next,  suppose $x_r=\max\limits_{1\leq i\leq n}{x_i}.$  Let $N_{\Gamma^{\prime}}(v_i)$ denote the  set of neighbours of  $v_i$ in $\Gamma^{\prime}$ and $d_{\Gamma^{\prime}}(v_i)$ denote the degree of $v_i$ in $\Gamma^{\prime}.$

\begin{claim}\label{duwei N-1}
$d_{\Gamma^{\prime}}(v_r)= n-1.$
\end{claim}
Otherwise, assume $d_{\Gamma^{\prime}}(v_r)\leq n-2$. Then
$$\lambda_1(\Gamma^{\prime})x_r=\sum_{v_j\in N_{\Gamma^{\prime}}(v_r)}\sigma^{\prime}(v_rv_j){x_j}\leq(n-2)x_r.$$ Thus, $\lambda_1(\Gamma^{\prime})\leq n-2,$
 a contradiction.

Without loss of generality, suppose $x_1\geq x_2$ and $x_3\geq x_4.$
\begin{claim}\label{r=3or4}
 $3\leq r\leq4.$
\end{claim}
Otherwise, assume $r\neq 3$ and $r\neq4.$ Recall that $x_r=\max\limits_{1\leq i\leq n}{x_i}$ and $d_{\Gamma^{\prime}}(v_r)=n-1.$  If $r>4,$  then there must exist an unbalanced $K_5$ induced by $\{v_1,v_2,v_3,v_4,v_r\},$ a contradiction.    If $r=1,$ then 
$$\begin{aligned}\lambda_1(\Gamma^{\prime})x_{r}&=\lambda_1(\Gamma^{\prime})x_{1}=-x_2+\sum_{3\leq s\leq n}x_s\\
&< 0+(n-2)x_r=(n-2)x_r.
\end{aligned}$$
Hence, $\lambda_1(\Gamma^{\prime})< n-2,$ a contradiction.  Thus, $r\neq1$ and $r\neq2.$ So Claim \ref{r=3or4} holds.

Recall that $x_3\geq x_4.$ Then by Claim \ref{r=3or4}, $x_3=\max\limits_{1\leq i\leq n}{x_i}.$
\begin{claim}\label{viadjacenttwo}
For any $5\leq i\leq n,$  at most two of $\{v_1,v_2,v_4\}$ are adjacent to   $v_i.$
\end{claim} 
Otherwise, assume that all of $\{v_1,v_2,v_4\}$ are adjacent to   $v_i.$ Note that $d_{\Gamma^{\prime}}(v_3)=n-1.$ Then there exists an unbalanced $K_5$ induced by $\{v_1,v_2,v_3,v_4,v_i\},$ a contradiction.

\begin{claim}\label{anyadjacent}
Any vertex of $\Gamma^{\prime}$ is adjacent to at least one  of $v_1$ and $v_2$.
\end{claim} 
Otherwise, without loss of generality, assume that $v_5$  is adjacent to neither $v_1$ nor $v_2.$ Then we can construct a new signed graph $\Gamma^{\ast}$  obtained from $\Gamma^{\prime}$ by adding a positive edge $v_1v_5.$  Clearly, $\Gamma^{\ast}$ is also unbalanced and $\mathcal{K}_5^-$-free. By Rayleigh principle, we obtain that
$$\begin{aligned}
\lambda_1(\Gamma^{\ast})-\lambda_1(\Gamma^{\prime})&\geq x^{\top}A(\Gamma^{\ast})x-x^{\top}A(\Gamma^{\prime})x\\
&=2x_1x_5>0.
\end{aligned} $$
Hence, $\lambda_1(\Gamma^{\ast})>\lambda_1(\Gamma^{\prime}),$ a contradiction.

Let $W_1=N_{\Gamma^{\prime}}(v_1)\backslash N_{\Gamma^{\prime}}(v_2)$  and $W_2=N_{\Gamma^{\prime}}(v_2)\backslash N_{\Gamma^{\prime}}(v_1).$ Recall that $x_1\geq x_2.$

\begin{claim}\label{W2=kong}
$W_2=\emptyset.$
\end{claim} 
Otherwise, assume $W_2\neq \emptyset.$ Without loss of generality, suppose $v_5\in W_2.$ Then we can construct a new signed graph $\Gamma^{\ast}$  obtained from $\Gamma^{\prime}$ by rotating the positive edge $v_2v_5$ to the non-edge position $v_1v_5.$  Clearly, $\Gamma^{\ast}$ is also unbalanced and $\mathcal{K}_5^-$-free. By Rayleigh principle, we obtain that
$$\begin{aligned}
\lambda_1(\Gamma^{\ast})-\lambda_1(\Gamma^{\prime})&\geq x^{\top}A(\Gamma^{\ast})x-x^{\top}A(\Gamma^{\prime})x\\
&=2x_5(x_1-x_2)\geq0.
\end{aligned} $$
If $\lambda_1(\Gamma^{\ast})=\lambda_1(\Gamma^{\prime})$, then  $x$ is also an eigenvector of $\Gamma^{\ast}$
corresponding to $\lambda_1(\Gamma^{\ast}).$ Based on the following equations,
$$\lambda_1(\Gamma^{\prime})x_{1}=\sum_{v_s\in N_{\Gamma^{\prime}}(v_1)}\sigma^{\prime}(v_sv_t)x_s$$ and
$$  \lambda_1(\Gamma^{\ast})x_{1}=\sum_{v_s\in N_{\Gamma^{\prime}}(v_1)}\sigma^{\prime}(v_sv_t)x_s+x_5,$$
we obtain that $x_5=0,$ which  contradicts  Claim \ref{1zero}.
Hence, $\lambda_1(\Gamma^{\ast})>\lambda_1(\Gamma^{\prime}),$ a contradiction.
\begin{claim}\label{W1adjacent all}
For any $v_p\in W_1,$  $v_p$ is adjacent all of $V(\Gamma^{\prime})\backslash\{v_p,v_2\}.$
\end{claim} 
Otherwise, assume  $v_q\in V(\Gamma^{\prime})\backslash\{v_p,v_2\}$ and $v_p\nsim v_q.$  Then we can construct a new signed graph $\Gamma^{\ast}$  obtained from $\Gamma^{\prime}$ by adding a positive edge $v_pv_q.$  Clearly, $\Gamma^{\ast}$ is also unbalanced and $\mathcal{K}_5^-$-free. By Rayleigh principle, we obtain that
$$\begin{aligned}
\lambda_1(\Gamma^{\ast})-\lambda_1(\Gamma^{\prime})&\geq x^{\top}A(\Gamma^{\ast})x-x^{\top}A(\Gamma^{\prime})x\\
&=2x_px_q>0.
\end{aligned} $$
Hence, $\lambda_1(\Gamma^{\ast})>\lambda_1(\Gamma^{\prime}),$ a contradiction.

Let $W=N_{\Gamma^{\prime}}(v_1)\cap N_{\Gamma^{\prime}}(v_2)\backslash\{v_3,v_4\}.$ 

\begin{claim}\label{independentset}
$W$ is an independent set.
\end{claim} 
Otherwise,  without loss of generality, assume $\{v_5,v_6\}\in W$ and $v_5\sim v_6.$ Note that $d_{\Gamma^{\prime}}(v_3)=n-1$ and  both $v_5$ and $v_6$ are adjacent to $v_1$ and $v_2.$  Then  there must exist an unbalanced $K_5$ induced by $\{v_1,v_2,v_3,v_5,v_6\},$ a contradiction.     
 
Let $w=|W|.$  Recall that $m=\frac{n(n-1)}{2}-h$ and    $h\leq n-3.$  
\begin{claim}\label{Wxiaoyu2}
$w\leq2.$
\end{claim} 
Otherwise, assume $w\geq3.$ Then by Claims \ref{viadjacenttwo} and  \ref{independentset},   $$\begin{aligned} h&\geq n-4+\frac{w(w-1)}{2}\geq n-4+3=n-1,
\end{aligned}$$
a contradiction.

\begin{claim}\label{w===0}
$w=0.$
\end{claim} 
Otherwise, by Claim \ref{Wxiaoyu2}, $w=1$ or $w=2.$  By Claim \ref{viadjacenttwo}, for any vertex $v_i\in W,$ $v_i\nsim v_4.$    If $w=1,$ then by Claims \ref{anyadjacent}, \ref{independentset}, \ref{W2=kong} and \ref{W1adjacent all}, $\Gamma^{\prime}=\Gamma_{1}$ (see Figure \ref{fig123}). However, by Lemma \ref{dayu n3}, $\lambda_1(\Gamma^{\prime})=\lambda_1(\Gamma_{1})<n-2,$  a contradiction. Thus, $w\neq 1.$  Similarly, $w\neq2.$ So Claim \ref{w===0} holds.

Above all,   $\Gamma^{\prime}=\Gamma_{3},$  which means $\Gamma$ is switching equivalent to $\Gamma_{3}$. This completes the proof.
\end{proof}

\end{document}